\documentclass[11pt]{amsart}

\usepackage[latin1]{inputenc}
\usepackage{amsmath,amsfonts,amssymb}
\usepackage{amscd}
\usepackage{graphicx}

\newcommand{\Ext}{\operatorname{Ext}}

\newcommand{\Hom}{\operatorname{Hom}}
\newcommand{\rad}{\operatorname{rad}}
\newcommand{\add}{\operatorname{add}}

\newcommand{\md}{\operatorname{mod}}

\newcommand{\cx}{\operatorname{cx}}
\newcommand{\id}{\operatorname{id}}

\newtheorem{theorem}{Theorem}[section]

\newtheorem{lemma}[theorem]{Lemma}

\begin{document}

\title[Maximal $\mbox{{\small $n$}}$-orthogonal 
modules for selfinjective algebras]
{Maximal $\mbox{{\large $n$}}$-orthogonal modules for
selfinjective algebras}

\author[Karin Erdmann, Thorsten Holm]
{Karin Erdmann \and Thorsten Holm}

\address{~~\newline
Karin Erdmann\newline
Mathematical Institute,
24-29 St. Giles,
Oxford OX1 3LB, U.K.}
\email{erdmann@maths.ox.ac.uk}

\address{~~\newline
Thorsten Holm\newline
Institut f\"ur Algebra und Geometrie,  
Otto-von-Guericke-Universit\"at Magdeburg, 
Postfach 4120, 39016 Magdeburg, Germany\newline
and\newline
Department of Pure Mathematics,
University of Leeds,
Leeds LS2 9JT, U.K.
}
\email{thorsten.holm@mathematik.uni-magdeburg.de}

\thanks{We gratefully acknowledge the support of the 
Mathematisches Forschungsinstitut Oberwolfach through 
a  Research in Pairs (RiP) project,
and also the support through a 
London Mathematical Society Scheme 4 grant.}

\bigskip

\begin{abstract}
Let $A$ be a selfinjective algebra.
We show that, for any $n\ge 1$, maximal $n$-orthogonal $A$-modules 
(in the sense of
Iyama) rarely  exist. More precisely, we prove that if $A$ admits a
maximal $n$-orthogonal module, then {\it all} $A$-modules are 
of complexity at most 1. 
\medskip

\end{abstract}

\maketitle

%%%%%%%%%%%%%%%%%%%%%%%%%%%%%%%%%%%%%%%%%%%%%%%%%%%%%%%%%%%%%%%%%

\section{Introduction}

Recently, O. Iyama introduced maximal $n$-orthogonal modules 
for finite-dimensional algebras, and developed an extensive
theory \cite{Iyama1}, \cite{Iyama2}. 
One aspect is a 'higher Auslander correspondence',
generalizing the famous one-one correspondence between algebras of
finite representation type and  Auslander algebras, that is, algebras
of global dimension at most 2 and dominant dimension
at least 2.

The existence of a maximal $n$-orthogonal module 
of an algebra $A$ has very striking consequences
for the homological properties of $A$ and its modules.
In particular it follows that then the representation dimension 
of $A$ (for background see \cite{Auslander}) is at most $n+2$. 
Of special interest is the case 
$n=1$. Maximal 1-orthogonal modules are known to exist 
for certain algebras of finite representation type, and also for
preprojective algebras \cite{GLS-rigid}. 
If $A$ has a 
maximal $1$-orthogonal module then
the representation dimension is at most 3. Using a result
of K. Igusa and G. Todorov \cite{IT}, this implies that the
famous finitistic dimension conjecture holds for $A$, that is, 
there is a finite bound on the projective dimensions of $A$-modules
of finite projective dimension.

If there are maximal $1$-orthogonal modules, then usually they
are not unique.
%there are many. 
However, Iyama showed that the endomorphism rings of any two maximal 
$1$-orthogonal modules of a fixed algebra are derived
equivalent (\cite{Iyama2}, 5.3.3). 
Moreover, he established a striking 'exchange rule':
taking  an indecomposable summand $X'$ of a maximal $1$-orthogonal
module, there is at most one indecomposable module $T$ not isomorphic
to $X'$ which can be substituted for $X'$ giving another maximal
1-orthogonal module (a proof can also be found in 
\cite{GLS-rigid}, 4.5).

Maximal
$1$-orthogonal modules are crucial for the work on cluster algebras of 
C. Gei{\ss}, B. Leclerc
and J. Schr\"oer \cite{GLS-rigid}. Cluster algebras were introduced
by Fomin and Zelevinsky in \cite{FZ} to study canonical bases of
quantum groups; a central feature is the introduction of an
'exchange graph'.  In the approach of \cite{GLS-rigid}, the exchange
property of maximal 1-orthogonal modules for preprojective
algebras describes the exchange graph for  the associated cluster algebra.
For details, see
\cite{GLS-rigid}.  

\smallskip

Because of these results, it would be very interesting to know 
how common maximal
1-orthogonal modules are. However we discovered 
that for selfinjective algebras they are very rare, and perhaps
occur only for the known cases of finite representation type and 
preprojective algebras. 
The aim of this note is to give a proof of this, and also
show that for any $n\ge 1$ maximal $n$-orthogonal modules are rare for 
selfinjective algebras.
This will also show that preprojective algebras play a very special role.
\medskip

We  recall the definition 
of a maximal $n$-orthogonal module for a finite-dimensional
algebra $A$, due to Iyama \cite{Iyama1}.
For  an $A$-module $X$, we  denote
by $\add(X)$ the full subcategory of the module category
${\rm mod}\,A$ whose objects are 
direct summands of direct sums of copies of $X$.
A (finitely generated)
$A$-module $X$ is called {\it maximal $n$-orthogonal} if 
for every $A$-module $M$ the following three conditions 
are equivalent
\begin{itemize}
\item[{(i)}] $\Ext^i_A(M,X)=0$ for all $1\le i\le n$.
\item[{(ii)}] $\Ext^i_A(X,M)=0$ for all $1\le i\le n$.
\item[{(iii)}] $M\in\add(X)$. 
\end{itemize}

For more details and some examples illustrating this concept, 
we refer to Section \ref{sec-max-n} below.
\smallskip

The following is the main result of this note, showing that 
only very few selfinjective algebras 
can possibly admit maximal $n$-orthogonal
modules. 

\begin{theorem} \label{thm-max1}
Let $A$ be a selfinjective algebra, and suppose that 
for some $n\ge 1$, there exists a maximal 
$n$-orthogonal $A$-module. Then  all 
$A$-modules have complexity at most 1. 
\end{theorem}

Recall that complexity of a module
measures the growth
of its  minimal projective resolution.
For selfinjective algebras, the most common modules
that have complexity $\leq 1$ are the $\Omega$-periodic
modules (here $\Omega M$ is the kernel of
a minimal projective cover of
the module $M$, and $M$ is $\Omega$-periodic if 
$\Omega^k(M)\cong M$ for some $k\geq 1$).

Let us point out that the existence of maximal $1$-orthogonal 
modules for finite-dimensional preprojective algebras is perfectly 
in line with our above results. 
In fact, for preprojective algebras of Dynkin type all modules
are $\Omega$-periodic, of period at most 6 (an 
unpublished result of C.\,M. Ringel and A. Schofield; 
for a proof see for instance \cite{ES1}, \cite{ES2},
or \cite{BBK}). 

\medskip

Algebras are in this paper assumed to be finite-dimensional algebras
over a field $K$. All modules are finitely
generated right modules, and \hskip-0.2cm $\mod A$ denotes the  category of
finitely generated $A$-modules.  

\section*{Acknowledgement}
We thank Jan Schr\"{o}er and 
{\O}yvind Solberg for helpful comments on some aspects
of this paper.

%%%%%%%%%%%%%%%%%%%%%%%%%%%%%%%%%%%%%%%%%%%%%%%%%%%%%%%%%%%
\section{Background and preliminaries}

\subsection{Homological algebra for selfinjective algebras}
\label{Ext}

Let $A$ be a finite-dimensional selfinjective algebra, so that 
projective modules and injective modules are the same. For a module
$M$, we have $\Omega M$, the kernel of a minimal projective cover, and
we also have $\Omega^{-1}M$, the cokernel of an injective hull.
Then $\Omega$ and $\Omega^{-1}$ induce mutually inverse equivalences
of the stable module category of $A$ (see for example \cite{ARS}, 
Chapter IV). Recall that the stable
module category $\underline{\md}\,A$ has the same objects as
${\rm mod}\,A$, and the morphisms  $\underline{\Hom}_A(M,N)$
are equivalence classes of module homomorphisms modulo those
factoring through a projective $A$-module. 
In particular we have for all $k\ge 1$ that
$$\Ext^k_A(M,N)\cong \underline{\rm Hom}_A(\Omega^kM, N)
\cong \underline{\rm Hom}_A(M, \Omega^{-k}N).
$$

\subsection{Maximal $n$-orthogonal modules} \label{sec-max-n}
For the convenience of the reader, we restate here Iyama's definition
of a maximal $n$-orthogonal module for a finite-dimensional
algebra, as already given in the 
introduction. 

\smallskip

An $A$-module $X$ is called {\it maximal $n$-orthogonal} if 
for every $A$-module $M$ the following three conditions 
are equivalent
\begin{itemize}
\item[{(i)}] $\Ext^i_A(M,X)=0$ for all $1\le i\le n$.
\item[{(ii)}] $\Ext^i_A(X,M)=0$ for all $1\le i\le n$.
\item[{(iii)}] $M\in\add(X)$. 
\end{itemize}

If $X$ is maximal $n$-orthogonal then
all projective indecomposable $A$-modules
and all injective indecomposable $A$-modules must be 
summands of $X$. Moreover,  $X$ does not have 
self-extensions, that is, $\Ext^i_A(X,X)=0$
for $i=1,\ldots,n$. 

We are interested in studying such modules when the
algebra $A$ is selfinjective but not semisimple. In that case,
$X$ must have at least one indecomposable summand which is not
projective (and injective). Namely, otherwise every indecomposable
$A$-module $M$ would have to be a summand of $X$ since 
${\rm Ext}^i_A(X, M)=0$ for $i=1,\ldots,n$ 
and then every indecomposable $A$-module would be projective and then even
simple, and $A$ would be semisimple. 

Furthermore, if $X$ is a maximal $n$-orthogonal 
module of a selfinjective algebra $A$ then so is $\Omega^t X\oplus A$
for any $t\in\mathbb{Z}$.

\subsubsection{Some examples}
Here are some easy explicit examples to illustrate the concept. 
\medskip

(1) First, let $A=K[T]/(T^t)$, a truncated polynomial ring.
Then $A$ has $t$ indecomposable modules (up to isomorphism), 
of dimensions $1,2,\ldots,t$. But all 
non-projective indecomposable $A$-modules have
self-extensions. Hence, there is no maximal $n$-orthogonal
$A$-module for any $n\geq 1$.
\medskip

(2) Let $Q$ be the following quiver 
\smallskip

%%%%%%%%%%%%%%%%%%%%%%

\begin{center}
\includegraphics[scale=0.5]{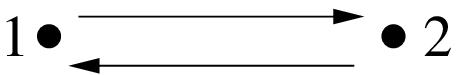}
\end{center}

%%%%%%%%%%%%%%%%%%%%%%
\smallskip

Set $A=KQ/\rad^2(KQ)$, a four-dimensional selfinjective
algebra. This algebra has precisely four
indecomposable modules, namely two simple modules $S_1$ and $S_2$, and
two indecomposable projectives  
$P_1$ and $P_2$ with $P_i$ the projective cover of $S_i$.
Let $X:= P_1\oplus P_2\oplus S_1$, then $X$ is a maximal $1$-orthogonal
$A$-module. In fact, $\Ext^1_A(X, S_2) \cong \Ext^1_A(S_1,S_2)\neq 0$ and 
$\Ext^1_A(S_2, X)\cong \Ext^1_A(S_2,S_1)\neq 0$. 

Then $X':= P_1\oplus P_2\oplus S_2$ also is 
a maximal $1$-orthogonal $A$-module, since $X' = \Omega X\oplus A$.  
\medskip

(3) We consider the quiver $Q$ as in (2),
and now let $A=KQ/\rad^3(KQ)$. Then the selfinjective algebra
$A$ has six indecomposable
modules, namely the simple modules $S_1, S_2$, their projective covers
$P_1$, $P_2$, and furthermore  
a 2-dimensional module $U_{1,2}$ with top $S_1$, and a 2-dimensional module
$U_{2,1}$ with top $S_2$. Note that $U_{1,2}\cong \Omega S_2$
and $U_{2,1}\cong \Omega S_1$. Suppose we have a
maximal $1$-orthogonal $A$-module $X$, then $X$ must have at least
one non-projective indecomposable summand. We may assume that
it has a simple summand (otherwise we replace $X$ by 
$\Omega^{-1}X\oplus A$).
At most one of the simples can be a summand of $X$ (since
$S_1$ and $S_2$ have a non-split extension). Suppose, say, 
$S_1$ is a summand
of $X$. Now, for any indecomposable non-projective module
$M\neq S_1$, we have
${\rm Ext}^1_A(M, S_1)\neq 0$ or ${\rm Ext}^1_A(S_1,M)\neq 0$. 
So $X$ can have no further non-projective
summands. On the other hand, 
$$\Ext^1_A(X, U_{1,2}) \cong 
\Ext^1_A(S_1,U_{1,2}) \cong \underline{\Hom}_A(U_{2,1},U_{1,2})=
0$$ a contradiction, since  $X$ is assumed to be  maximal
1-orthogonal. 

\medskip

(4) For any natural number $n$ there exists a selfinjective
algebra $A$ with an $n$-orthogonal module, as the
following example shows. 

For any $n\ge 1$, let $Q$ be the cyclic (oriented) quiver 
with $m:=2n+2$ vertices. Then consider the selfinjective algebra
$A=KQ/\rad^2(KQ)$. Note that the indecomposable $A$-modules
are the projectives $P_0,P_1,\ldots,P_{m-1}$, and the simple
modules $S_0,S_1,\ldots,S_{m-1}$. We label the simple modules so that
$P_i$ has socle $S_{i+1}$, with indices taken modulo $m$. 
Then it is straightforward to check that the module
$$X:=P_0\oplus P_1\oplus\ldots\oplus P_{m-1}\oplus S_0\oplus S_{n+1}$$
is a maximal $n$-orthogonal $A$-module. 
\bigskip

%The above easy examples should indicate that in general 
In general, it is not at all easy to decide whether or not
maximal $n$-orthogonal modules exist. 
In \cite{Iyama2}, O. Iyama discusses the case $n=1$ for
selfinjective algebras of finite representation type,
and gives a combinatorial reformulation in terms
of certain triangulations of regular $m$-gons. 
For preprojective algebras,  
the existence of maximal $1$-orthogonal modules is proved
in \cite{GLS-rigid}.

\subsection{Auslander-Reiten formula \cite{AR}} 
\label{AR-formula}
Let $A$ be any 
finite-dimensional algebra, and let $\tau=D\,Tr$ be the 
Auslander-Reiten translation. When $A$ is selfinjective, we have
$\tau \cong \Omega^2 \nu$ where $\nu$ is a Nakayama automorphism
of $A$, see \cite{ARS} Chapter IV, 3.7. 
Then for any $A$-modules 
$M, N$ we have
$$D\Ext^1_A(M,N) \cong \underline{\Hom}_A(\tau^{-1}N, M).$$
We will use freely that $\Omega$ and $\nu$ commute for a selfinjective
algebra $A$.

%%%%%%%%%%%%%%%%%%%%%%%%%%%%%%%%%%%%%%%%%%%%%%%%%%%%%%%%%%%%
\section{Periodicity of $n$-orthogonal modules}

The following theorem is the first crucial step in proving
our main result. 

\begin{theorem} \label{thm-n-orth}
Let $A$ be a selfinjective algebra, and, for some $n\ge 1$, let
$X$ be a maximal $n$-orthogonal $A$-module. If $Y$ is a direct 
summand of $X$ then so is $\Omega^{n+2}\nu Y$. 
Hence 
every non-projective indecomposable
summand of $X$ is $\Omega^{n+2}\nu$-periodic. 
\end{theorem}

Before embarking on the proof, was make an easy but
useful observation.

\begin{lemma} \label{lem-Ext}
Let $A$ be a selfinjective algebra. For any 
$A$-module $M$ and any $i\ge 1$ we have an isomorphism of vector spaces
$$\Ext^i_A(M,N) \cong \Ext^1_A(N,\Omega^{i+2}\nu M).$$
\end{lemma}

\proof Using the Auslander-Reiten formula \ref{AR-formula}
and the formula \ref{Ext} we get 
\begin{eqnarray*}
\Ext^{i}_A(M, N) & \cong  & \underline{\Hom}_A(\Omega^{i}M, N)
 \cong  \underline{\Hom}_A(\Omega M,\Omega^{-i+1}N) \\
  & \cong & \Ext^1_A(M,\Omega^{-i+1}N) \cong
     \underline{\Hom}_A(\tau^{-1}\Omega^{-i+1}N,M) \\
  & \cong & \underline{\Hom}_A(\Omega N,\tau\Omega^iM)
 \cong \Ext^1_A(N,\tau\Omega^iM) \\
 & \cong & \Ext^1_A(N,\Omega^{i+2}\nu M)
\end{eqnarray*}
where for the last isomorphism we use the fact that,  for $A$ selfinjective
algebra, one has $\tau=\Omega^2\nu$. 
\qed
\medskip

Now we are in the position to complete the proof of Theorem
\ref{thm-n-orth}. 
\medskip

{\it Proof of Theorem \ref{thm-n-orth}.}
Let $X$ be a maximal $n$-orthogonal module for the selfinjective
algebra $A$. We consider the $A$-module $\Omega^{n+2}\nu X$.
For any $i$ such that $0< i\le n$ we obtain
\begin{eqnarray*}
\Ext^i_A(X,\Omega^{n+2}\nu X) & \cong & 
  \underline{\Hom}_A(\Omega^{i}X,\Omega^{n+2}\nu X) 
 \cong  \underline{\Hom}_A(\Omega X,\Omega^{-i+n+3}\nu X) \\
 & \cong & \Ext^1_A(X,\Omega^{-i+n+3}\nu X) 
 \cong  \Ext^{-i+n+1}_A(X,X) \\
\end{eqnarray*} 
where the last isomorphism comes from Lemma \ref{lem-Ext}. Note that the 
superscripts $-i+n+1$ run through the set $\{1,2,\ldots,n\}$. Since
$X$ is maximal $n$-orthogonal, we conclude that 
$\Ext^{-i+n+1}_A(X,X)=0$ for all $i=1,\ldots,n$, which means that
$\Ext^i_A(X,\Omega^{n+2}\nu X) = 0$ for $i=1,\ldots,n$. 
Using again that $X$ is maximal $n$-orthogonal we deduce
that $\Omega^{n+2}\nu X \in \add(X)$.

In particular, if $Y$ is an indecomposable direct summand of $X$,
then also $\Omega^{n+2}\nu Y$ is a direct summand of $X$. 
\smallskip
 
This means that $\Omega^{n+2}\nu$ permutes the indecomposable
non-projective summands of $X$
(recall that $\Omega^{n+2}\nu$ induces a permutation on the 
set of non-projective indecomposable $A$-modules). 
But $X$ has by definition only finitely many indecomposable 
summands. Hence some power of $\Omega^{n+2}\nu$ is the identity
permutation on the non-projective summands of $X$, that is, 
$X$ is $\Omega^{n+2}\nu$-periodic. 
\qed

%%%%%%%%%%%%%%%%%%%%%%%%%%%%%%%%%%%%%%%%%%%%%%%%%%%%%%%%
\section{Complexity at most 1}

\subsection{Complexity} Let $A$ be a finite-dimensional algebra. 
For any $A$-module $M$, let
$$\ldots \to P_2\to P_1\to P_0\to M\to 0$$
be a minimal projective resolution. The complexity of $M$
measures the rate of growth of the terms
of such a resolution. More precisely, the {\em complexity of
M} is defined as
$$\cx(M):=\inf \{b\in\mathbb{N}_0\mid \exists\, c>0: 
\dim P_n\le c\, n^{b-1}\mbox{ for all $n$}\},$$
if it exists, otherwise $\cx(M)=\infty$. 
Note that $\cx(M)=0$ precisely for modules $M$ having 
finite projective dimension. Moreover, we have $\cx(M)\le 1$ 
if and only if the dimensions of the $P_n$'s are bounded.
Clearly, if $M$ is $\Omega$-periodic then $\cx(M)=1$. 
The converse is not true in general, see \cite{LS} for a
counterexample.
\medskip

The following well-known result will be useful later. 

\begin{lemma} \label{lem-ses-cx}
Let $A$ be a finite-dimensional algebra. 
\begin{enumerate}
\item[{(a)}] 
Suppose 
$$0\to L\to M\to N\to 0$$ is a short exact sequence
of $A$-modules. Then, if two of the modules in the sequence
have complexity $\le 1$, then so does the third. 
\item[{(b)}] Let $\alpha$ be an automorphism of the algebra 
$A$, and for any $A$-module $M$ let $M_{\alpha}$ denote 
the $A$-module with twisted action $m\cdot a:=m\alpha(a)$. 
Then $M$ and $M_{\alpha}$ have the same complexity.  
\end{enumerate}
\end{lemma}

To prove Theorem 1.1, we will use 
the following
standard construction which we recall here for the convenience of 
the reader.

\begin{lemma} \label{lem-universal}
(Universal extension)
Let $A$ be a finite-dimensional algebra, and let $X$ be an 
$A$-module with $\Ext^1_A(X,X)=0$. Moreover, let $V$ be an
$A$-module such that 
$n:=\dim \Ext^1_A(X,V)>0$.
Then there exists a short exact sequence
$$ 0 \to V \to U \to X^n\to 0$$
for which $\Ext^1_A(X,U)=0$.
\end{lemma}

\proof (\cite{Bongartz}, Lemma 2.1)
We choose a basis of $\Ext^1_A(X,V)$, say $e_1,\ldots,e_n$. 
Then we construct the short exact sequence
$$0\to V\to U\to X^n\to 0$$
such that the pullback under the $i^{th}$ canonical injection
$X\to X^n$ is $e_i$. 
Upon applying $\Hom_A(X,-)$ we get a long exact sequence
$$\ldots \to \Hom_A(X,X^n) \stackrel{\delta}{\twoheadrightarrow} 
\Ext^1_A(X,V) \to \Ext^1_A(X,U) \to \underbrace{\Ext^1_A(X,X^n)}_0 
\to \ldots
$$
Note that the map $\delta$ is surjective by construction. 
Hence, it follows that $\Ext^1_A(X,U)=0$, as desired.
\qed

\medskip

{\it Proof of Theorem 1.1.} 
Any module which is  $\Omega^{n+2}\nu$-periodic 
has complexity $\leq 1$. So we assume now that $A$ has an 
indecomposable module $V$ which is not 
$\Omega^{n+2}\nu$-periodic. 

Let $V=U_0$.
We construct inductively modules
$U_1, U_2, \ldots, U_n$ and short exact sequences
$$0\to U_{i-1} \to U_i\to (\Omega^{n-i}X)^{r_i}\to 0
\leqno{(\zeta_i)}$$
for $1\leq i \leq n$, such that
${\rm Ext}^1_A(\Omega^jX, U_i)=0$ for $n-i\leq j\leq n-1$. 

\medskip

(a) We first construct $U_1$. 
If ${\rm Ext}^1_A(\Omega^{n-1}X, V)=0$ then take 
$U_1 = V\oplus \Omega^{n-1}X$ and $r_1=1$. 
Otherwise, construct the universal extension (see 
\ref{lem-universal})
$$0\to V\to U_1\to (\Omega^{n-1}X)^{r_1}\to 0.
$$
Then ${\rm Ext}^1_A(\Omega^{n-1}X, U_1)=0$.

\medskip

(b) For the inductive step, suppose $U_1, \ldots, U_{i-1}$ have been
constructed. 

If ${\rm Ext}^1_A(\Omega^{n-i}, U_{i-1})=0$ then take
$U_i = U_{i-1}\oplus \Omega^{n-i}X$ and $r_i=1$. 
Otherwise, we construct the universal extension
$$0\to U_{i-1}\to U_i\to (\Omega^{n-i}X)^{r_i}\to 0.
$$
Then by construction we have 
${\rm Ext}^1_A(\Omega^{n-i}X, U_i)=0$. 
Furthermore, for $n-i<j$ we have by the inductive hypothesis that
${\rm Ext}^1_A(\Omega^jX, U_{i-1})=0$.

Since $X$ is maximal $n$-orthogonal we have
for all $k=0,1,\ldots,n-1$ that
$$\Ext^1_A(\Omega^k X, X)= \Ext^{k+1}_A(X,X)=0.$$

In particular, in our situation we then 
know that 
$${\rm Ext}^1_A(\Omega^jX, \Omega^{n-i}X) = 
\Ext^1_A(\Omega^{j-(n-i)}X,X) = 0.$$
By considering the long exact sequence to the previous
universal extension, one concludes
that ${\rm Ext}^1_A(\Omega^jX, U_i)=0$,
thus completing the inductive step.

\medskip

The module $U_n$ satisfies 
$$
\Ext^{j+1}_A(X,U_n) = {\rm Ext}^1_A(\Omega^jX, U_n)=0$$
for 
$0\leq j\leq n-1$. Since $X$ is maximal $n$-orthogonal 
it follows that 
$U_n$ belongs to ${\rm add}(X)$. In particular, $U_n$ is 
$\Omega^{n+2}\nu$-periodic, by 
Theorem \ref{thm-n-orth},
and then has complexity $\leq 1$, by Lemma 
\ref{lem-ses-cx}.
Now we use downward induction. In the extension
$(\zeta_n)$, the last two terms
have complexity $\leq 1$ and hence so does the first term, that is,
$U_{n-1}$, by Lemma \ref{lem-ses-cx}. For the inductive step, suppose
$U_j$ has complexity $\leq 1$, then the last two terms in 
the sequence $(\zeta_{j})$ have complexity $\leq 1$ and hence so does
$U_{j-1}$, again by \ref{lem-ses-cx}. 
The last step shows that
$V$ has complexity $\leq 1$.
\qed

\medskip

\section{Concluding remarks and open questions}

\subsection{} 
We have proved that if $A$ is selfinjective and $A$ has a maximal
$n$-orthogonal module then all $A$-modules have complexity
$\leq 1$. 
For algebras of finite representation 
type, this is no restriction. On the other hand, it is a very strong
restriction in general.
One would like to know which algebras of infinite 
representation type have this
property, and whether any such algebra has a maximal $n$-orthogonal
module. 
\smallskip

Algebras for which every module has $\Omega^3\nu$-period $\leq 2$ were
classified in \cite{BES}. The list consists of the preprojective
algebras of Dynkin type, then one series of algebras 
denoted by $P(\mathbb{L}_n)$ (where $n\ge 2$)
which have precisely one simple module with self-extensions, and 
otherwise certain deformations of these algebras. By \cite{GLS-rigid}, 
preprojective
algebras do have maximal 1-orthogonal modules. On the other
hand, they do not have maximal $n$-orthogonal modules
for $n\geq 2$ since $\Ext^2(M,M)\neq 0$ for all
non-projective modules. It is also easy to
see that $P(\mathbb{L}_2)$ (which is of finite type) does not
have  a maximal $1$-orthogonal module. 
We do not know whether or not $P(\mathbb{L}_n)$ for $n\geq 3$, or the 
deformations of preprojective
algebras have maximal 1-orthogonal modules.
\smallskip

Tame selfinjective algebras for which
all modules are periodic can be found in 
\cite{BS1}, \cite{BS2}. Moreover, there are as well
the algebras of quaternion type in \cite{Erdmann}.

\subsection{}
For algebras $A$ of quaternion type, we can 
show that they can not have a maximal
$n$-orthogonal module, for any $n\neq 2$. 
Recall that by definition $A$ is symmetric (hence
$\tau\cong \Omega^2$ and $\nu\cong \id$) and for any $A$-module $M$ one has
$\Omega^4M\cong M$ (see \cite{Erdmann}). 
First, for any non-projective indecomposable 
$A$-module $M$ we get from Lemma \ref{lem-Ext}
\begin{eqnarray*}
\Ext^3_A(M,M) \cong  \Ext^1_A(M,\Omega^5M) \cong
 \Ext^1_A(M,\Omega M) \cong \underline{\Hom}_A(\Omega M,
\Omega M) \neq 0.
\end{eqnarray*}
In particular, $A$ cannot have a maximal $n$-orthogonal
module where $n\ge 3$.

Secondly, suppose $X$ is a maximal $1$-orthogonal $A$-module.
Then by Lemma \ref{lem-Ext} we have
$$\Ext^1_A(X,\Omega^{-1}X) \cong \Ext^1_A(X,X) =0.$$
Hence, for any non-projective indecomposable summand 
$Y$ of $X$, we also have $\Omega^{-1}Y\in \add(X)$.
But on the other hand,
$$\Ext^1_A(\Omega^{-1}X,X) \cong 
\underline{\Hom}_A(X,X) \neq 0,$$
which contradicts the maximal $1$-orthogonality of $X$. 
\medskip

\subsection{} We do not know any algebra of infinite type for which all
modules have complexity $\leq 1$ but which has modules which are
not $\Omega^{n+2}\nu$-periodic. 
Such algebra, if it exists, 
would have very unusual homological properties.
For example if for such algebra, the Nakayama automorphism
$\nu$ has finite order,
then the finite generation properties Fg1 and Fg2 in \cite{EHSST} must
fail, since those imply that $\Omega$-periodicity is the same as 
having complexity one.

\subsection{} Maximal $n$-orthogonal modules for an algebra $A$
have all projective indecomposable and all injective
indecomposable $A$-modules as direct summands. Since
$\Ext^i_A(X,X)=0$ for $i=1,\ldots,n$ for any 
maximal $n$-orthogonal module, there cannot be any non-split 
extensions between injective and projective $A$-modules. 
Therefore, when searching for maximal $n$-orthogonal modules
it seems very natural to consider selfinjective algebras,
as we did in this paper. However, there are 
non-selfinjective algebras having maximal $n$-orthogonal 
modules. The  examples we know are all of
finite type, and hence its modules also are periodic. It would be
interesting to know
whether there can exist maximal $n$-orthogonal
modules for non-selfinjective algebras $A$ for which not
all $A$-modules are of complexity at most 1.

%\medskip

%%%%%%%%%%%%%%%%%%%%%%%%%%%%%%%%%%%%%%%%%%%%%%%%%%%%%%%%%%%%%%

\end{document}